\documentclass[12pt,reqno]{amsart}
\usepackage{amssymb}
\usepackage{amscd}
\usepackage{amsxtra}
\usepackage{amsmath}
\usepackage[mathscr]{eucal}
\usepackage{tikz-cd, enumitem}
\usepackage{xurl}
\usepackage{hyperref}

\theoremstyle{plain}
\newtheorem{thm}{Theorem}

\theoremstyle{definition}

\newtheorem{ques}{Question}

\theoremstyle{remark}

\newcommand{\bbH}{\mathbb{H}}
\newcommand{\Z}{\mathbb{Z}}

\newcommand{\PGL}{\mathsf{PGL}}

\newcommand{\PSL}{\mathsf{PSL}}

\begin{document}
\title[]{An improper surface group action}
\author{Nic Brody}
\date{\today}
\maketitle

\begin{abstract}
	In 1999, Long and Reid proposed a proper action of a surface group on a product of trees. In this note, we show that the action is not proper.
\end{abstract}

The following question \cite{BFMV25, Button19, FLSS18, Schwartz23} relates two of the most fundamental objects of consideration in low-dimensional topology and geometric group theory. 

\begin{ques}
	Can a surface group act properly on a locally compact product of trees?
\end{ques}

A folklore candidate action, first considered by Long and Reid in the 1990s, attempted to produce an explicit example of this. Let $\Sigma$ denote the orbifold which is a torus with a single cone point of order 2. Its fundamental group is given by the presentation $\pi_1(\Sigma)=\langle a,b\mid [a,b]^2=1\rangle$. Magnus \cite{magnus} considers the $\PGL_2$ character variety of this group. An element of $\PGL_2$ has order 2 if and only if it has trace 0, and so we seek two matrices whose commutator has trace zero to obtain a representation of this group. Supposing the first matrix is $a=\begin{pmatrix} t & 0 \\ 0 & 1 \end{pmatrix}$, the second matrix $b$ must satisfy $2b_{11}b_{22}=(t+t^{-1})b_{12}b_{21}$. We obtain 

\[ \left \langle \begin{pmatrix} t & 0 \\ 0 & 1\end{pmatrix}, \begin{pmatrix} 1+t^2 & 2 \\ t & 1 \end{pmatrix}\right \rangle \leq \PGL_2\left(\Z[t]\left[\frac{1}{t},\frac{1}{t-1}\right]\right).\]

Any choice of real parameter $t>1$ determines a faithful representation of this group, because one can construct an explicit hyperbolic rectangle with this geometric structure. We now specialize to $t=9$, and obtain the Long-Reid group 

\[ \Gamma=\left \langle \frac{1}{3}\begin{pmatrix} 9 & 0 \\ 0 & 1\end{pmatrix}, \frac{1}{8}\begin{pmatrix} 82 & 2 \\ 9 & 1 \end{pmatrix}\right \rangle \leq \PSL_2\left(\Z\left[\frac{1}{6}\right]\right).\]

By \cite{BHC}, $\PGL_2(\Z[1/6])$ is a lattice in $\bbH^2 \times T_3 \times T_4$, where $T_k$ denotes the $k$-regular tree. The stabilizer of a (particular) vertex in $T_3\times T_4$ is $\PGL_2(\Z)$, and the action of $\Gamma$ on $T_3 \times T_4$ is proper if and only if $\Gamma\cap \PGL_2(\Z)$ is finite.

\begin{thm}
	The action of $\Gamma$ on $T_3\times T_4$ is not proper.
	\begin{proof}
		The length 82 word 
		\[
		\begin{aligned}
				& a^2ba^{-2}b^{-1}aba^2b^{-1}a^{-1}b^{-1}a^2ba^{-1}ba^2b^{-1}a^{-1}b^{-1}abab^{-1}aba^{-1}ba^2b^{-1}a^3ba^{-2}b^{-1} \\
				& ab^{-1}a^{-2}b^2ab^{-1}a^{-2}ba^{-1}b^{-1}ab^{-1}a^{-2}bab^{-1}a^2ba^{-1}bab^{-1}aba^{-1}ba^2b^{-1}a^3ba^{-2}b^{-1}
		\end{aligned}
		\]
		
		evaluates to the infinite-order integer matrix \[\begin{pmatrix}-646279884109511971664607 & 6162511442411222450262052 \\ -4193268331567764626734 & 39984323680432243295081\end{pmatrix}.\]
	\end{proof}	
\end{thm}

\bibliographystyle{amsplain}
\bibliography{LongReidGroup.bib}
\nocite{*}


\end{document}